\documentclass[a4paper,11pt,reqno]{amsart}

\usepackage{a4wide}
\usepackage{amsthm}                   
\usepackage{amsmath}                  
\usepackage{amsfonts}                 
\usepackage{amssymb}                  
\usepackage{amsbsy}                   
\usepackage{mathrsfs}                 
\usepackage{mathbbol}                 
\usepackage{graphicx}                 
\usepackage{fontenc}                  
\usepackage{subfigure}                   
\usepackage[isolatin]{inputenc}       
\usepackage{bbm}

\newtheorem*{zeroinfty}{$0$-$\infty$-Law of Stochastic Geometry}
\newtheorem{proposition}{Proposition}

\newcommand{\colwidth}{0.7\textwidth}
\newcommand{\smallwidth}{0.3\textwidth}
\newcommand{\halfwidth}{0.4\textwidth}

\newcommand{\e}{\ensuremath{\mathrm{e}}}      
\newcommand{\D}{\ensuremath{\mathscr{D}}}     
\newcommand{\dr}{\ensuremath{\mathscr{D}_R}}  
\newcommand{\dv}{\ensuremath{\mathscr{D}_V}}  
\newcommand{\ds}{\ensuremath{\mathscr{D}_S}}  
\newcommand{\R}{\ensuremath{\mathbb{R}}}      
\newcommand{\Rd}{\ensuremath{\R^d}}           
\newcommand{\mdrd}{\ensuremath{{\mathscr{M}^{\bullet}\big(\mathbb{R}^d\big)}}}
\newcommand{\mdrn}{\ensuremath{{\mathscr{M}^{\bullet}\big(\mathbb{R}^n\big)}}}
\newcommand{\mdrt}{\ensuremath{{\mathscr{M}^{\bullet}\big(\mathbb{R}^2\big)}}}
\newcommand{\mdg}{\ensuremath{{\mathscr{M}^{\bullet}(\mathfrak{X}_d)}}}

\newcommand{\x}{\ensuremath{X}}

\DeclareMathOperator{\card}{card}             
\DeclareMathOperator{\supp}{supp}             
\DeclareMathOperator{\vol}{vol}             

\usepackage[pdfpagelabels,a4paper,colorlinks]{hyperref}
\hypersetup{%
   colorlinks=false,        
   linkcolor=black,
   backref=false,
   pagebackref=false,
   pdfpagemode=UseThumbs, 
   bookmarksopen=true,     
   bookmarksnumbered=true, 
   pdfstartpage={1},       
   pdftitle={Random Cluster Tessellations},
   pdfsubject={Construction of cluster tessellations by means of
     cluster properties and point processes},
   pdfauthor={Kai Matzutt},%
   pdfkeywords={Delone, Voronoi,
    random tesselations, point processes, random tiling, poisson process},
   pdfcreator={TeX},%
   pdfproducer={LaTeX with hyperref and thumbpdf}
   plainpages=false,
   hypertexnames=true
}

\begin{document}

\title{Random Cluster Tessellations}
\author{Kai Matzutt}
\address{Fakult\"at f\"ur Mathematik\\
Universit\"at Bielefeld\\
Pf.~100131\\
33501 Bielefeld\\
Germany
}
\date{\today}
\curraddr{}
\email{kai@math.uni-bielefeld.de}
\urladdr{http://www.math.uni-bielefeld.de/\raisebox{-0.8ex}{\~{ }}kai}
\dedicatory{}
\begin{abstract}
  This article describes, in elementary terms, a
  generic approach to produce 
  discrete random tilings and similar random structures by using point
  process theory. The standard
  Voronoi and Delone tilings can be constructed in
  this way. For this purpose, convex polytopes are replaced by their
  vertex sets. Three explicit constructions are given to illustrate the
  concept. 
\end{abstract}

\maketitle

\section{Introduction}

Apart from symmetry and long-range order, also randomness is needed to
provide appropriate structure models in physics or chemistry. 
A well-known example for an application of \emph{random points}
respectively \emph{point processes} is an
ideal gas, which, at each instance of time, may be 
described by a \emph{Poisson point process} (cf.~\cite{heer72}
p.~449f).

\begin{figure}[hbt]
    \centering
    \includegraphics[width=\halfwidth]{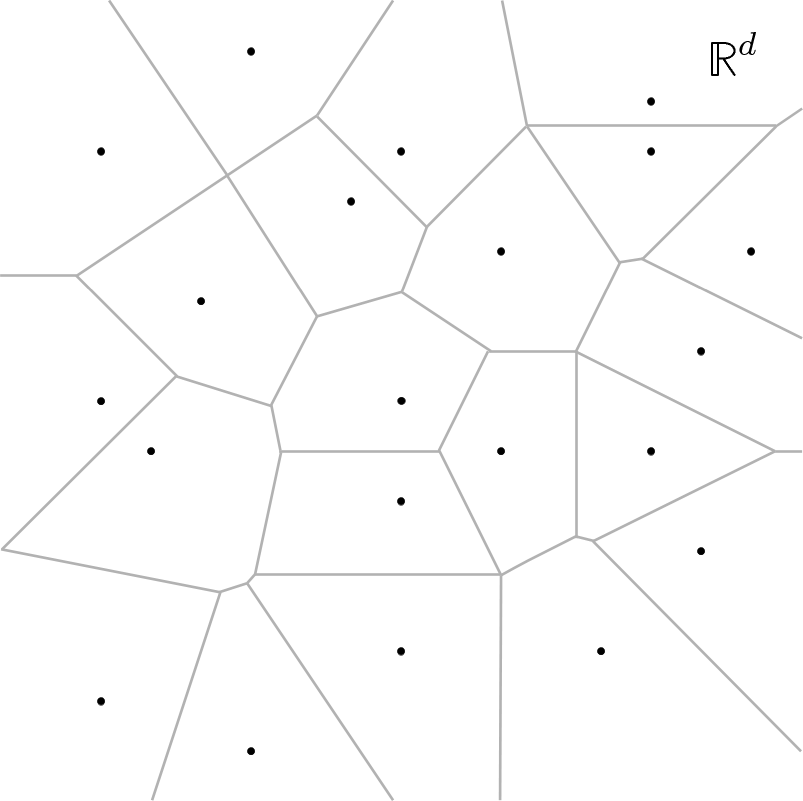}
    \caption{A Voronoi tessellation: The (random) points generate 
      their `neighbourhoods'}
    \label{fig:voronoi}
\end{figure}

There are also some well-known random tilings with applications in
crystallography and material sciences, such as the Poisson Voronoi and
the Poisson Delone tessellations (\cite{skm95,okabe92}). A more
general description of random tilings in the context of quasicrystals
can be found in \cite{rhhb98,henley99}.  Recent work
of Gummelt \cite{gummelt04,gummelt06} is also concerned with
establishing a connection between randomness and
quasicrystalline structures. 
Most approaches start from randomly generated points in Euclidean space
$\Rd$ and then distill the information about the tiles, often convex
polytopes, out of the positions of the points. (See Figure
\ref{fig:voronoi} for an illustration of this widely used concept.)  

\begin{figure}[hbt]
  \centering
   \includegraphics[width=\halfwidth]{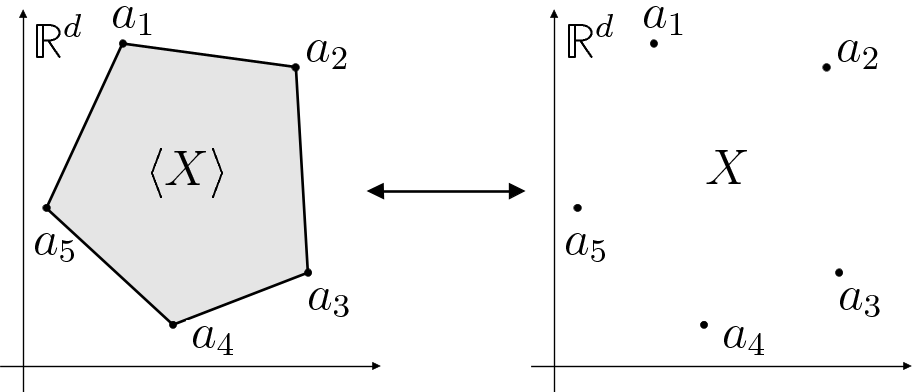}
   \caption{We identify a convex polytope with its vertex set. Here,
     $\langle \x\rangle$ denotes the convex hull of the set 
     $\x=\{a_1,\ldots,a_5\}\subseteq\Rd$}
   \label{fig:vertex-set}
\end{figure}

The theory presented in this article is partly based on recent work of Zessin
\cite{zessin05}, where random tiles are also extracted from random
point conf\/igurations. However, this approach replaces the concept of
convex polytopes by
discrete, f\/inite subsets of $\Rd$, called
\emph{clusters}. Nevertheless, tilings 
of convex polytopes can be embedded into this theory by
identifying the vertex set of a convex polytope with the polytope
itself (see Figure \ref{fig:vertex-set}). This point of view might
describe the underlying structure of the atoms and molecules of
certain materials more accurately. 
In this article the desired local and global properties of the tiles,
in our case the clusters, enter by means of so called \emph{cluster
  properties}. In particular, the f\/irst two examples given, are to illustrate
this concept.  
Also, although the examples given in this article are tilings, the
presented theory is not restricted to these.

While the typical point processes like the \emph{Poisson
  point process} with constant intensity might be 'too random' to
describe condensed matter, this paper gives an
explicit example construction
to go over to a random tiling 'close' to a quasicrystalline
one.


\section{Random Points}
Point processes can be constructed in various spaces. Here, we 
restrict ourselves to random points in $\Rd$. See \cite{kmm78} for a
more detailed and more general description.
To speak about randomness and probabilities, we f\/irst need to
identify the objects which should be realized randomly. Then we need
a notion of which events can be computed or, more precisely, can be
given a probability.

Technically, the space for the point
conf\/igurations is the set of the locally f\/inite simple counting
measures in $\Rd$, denoted by $\mdrd$. All possible events are
subsumed under the $\sigma$-algebra generated by the counting
functions. (A very readable introduction to probability theory and thereby
an explanation for the need of terms like measurability and
$\sigma$-algebras is given in \cite{georgii04}.) A probability measure 
on $\mdrd$ equipped with this $\sigma$-algebra is called
\emph{(simple) point process}.  
 
In essence, we can say that the point conf\/igurations we want
to consider have only f\/initely many points in every bounded subset of
$\Rd$, and that the typical events are of the form
\begin{equation}\label{eq:typ-event}
  \left\{\eta\in\mdrd:\,\eta\text{ has }k\text{
       points in }B\right\}\,,
\end{equation}
where $B$ is some bounded subset of $\Rd$ and $k$ is a non-negative
integer. 
A point process may be viewed as a mechanism to randomly generate
point conf\/igurations obeying a given probability law.
It is comprehensible that the probabilities of those events describe
the properties of random point sets in great detail since the bounded
sets $B$ in \eqref{eq:typ-event} can be chosen arbitrarily
small.
For calculations, it is very convenient to express a point
conf\/iguration  $\eta$ as a sum of Dirac measures,
\begin{equation*}
  \eta=\sum_{i=1}^{\infty}\delta_{a_i}.
\end{equation*}
Since $\eta$ is assumed to be locally f\/inite, one can always f\/ind a
suitable sequence $a_i$, $i=1,2,\ldots$ of points in $\Rd$, e.g., by
collecting the points in centred balls of increasing radius and giving
them consecutive labels.

An interesting class of point processes are the stationary ones, where
the probabilities of the events 
\eqref{eq:typ-event} are
invariant under translations of the bounded sets $B$. 
The best explored class of point processes is the class of
\emph{Poisson point processes} and among them the stationary ones in particular:
Let $\lambda\in\R^+$ and let $\vol 
(B)$ denote the volume of a Borel set $B$. The 
\emph{Poisson point process with intensity $\lambda$}, $P_\lambda$, 
assigns to our typical events the probabilities
\begin{equation}
  \label{eq:1}
  P_\lambda\left(\left\{\eta\in\mdrd\text{ has }k\text{
        points in }B\right\}\right)
  =\frac{\left(\lambda\vol
    (B)\right)^k}{k!}\exp\left(-\lambda\vol (B)\right)\,.
\end{equation}
The expected number of points in a unit cube thus equals
$\lambda$. It is clear that $P_\lambda$ is stationary
because the probabilities only depend on the translation invariant
volume. 

In the general case of Poisson point processes, the volume in
\eqref{eq:1} is replaced
by an arbitrary locally f\/inite measure $\rho $ on $\Rd$. This
results in the Poisson point process with \emph{intensity measure} $\rho $,
denoted by $P_\rho $. The process is
stationary as long as the intensity measure is translation invariant,
which in the case of $\Rd$ just means that $P_\rho=P_\lambda$ for some
given $\lambda>0$.
For special intensities, it might happen, with positive
probability, that point conf\/igurations have more than one point in one
place. Such situations are excluded if the intensity measure $\rho$
has no pure point part, i.e., if $\rho\left(\{a\}\right)=0$ for all
$a\in\Rd$. Such point processes are called \emph{simple}. 
While $P_\lambda$ describes some ideal gas, one might interpret
$P_\rho$ as gas of non-interacting molecules in a certain physical 
potential. In contrast to general \emph{Gibbs measures}
(cf.~\cite{georgii88}), the particles in the
Poissonian case are always non-interacting.


\section{Clusters and Cluster Properties}\label{sec:clusters}
Similar to the above mentioned models for an ideal gas, randomness
enters the approach of this paper by means of point point processes,
as described in the previous section, where there are a lot of
well-known constructions and simulations \cite{schneider-weil00, skm95}. But the information one gets out
of the typical construction rules are of a more global nature, like
distributions. To describe the local properties of the modelled
objects and still not to lose the randomness of the point processes, we
need a proper concept, which will be described in this section. 

As mentioned before, a cluster from our point of view is a
f\/inite subset of $\Rd$. Let
\begin{equation*}
  \mathfrak{X}_d:=\left\{\x\subset \Rd\,\big\vert\,\card (\x)<+\infty\right\}
\end{equation*}
be the space of clusters in $\Rd$ (here $\card (\x)$ denotes the
cardinality of a set $\x$). Typical events in this space are
constructed analogously to $\mdrd$.
The method to consider (random) collections of clusters is inspired
by the theory of \emph{random sets} by Matheron \cite{matheron75}. 

To attach clusters $\x\in\mathfrak{X}_d$ to point conf\/igurations
$\eta\in\mdrd$, we use the concept of \emph{cluster properties}.
A cluster property $\mathscr{D}$ is a measurable (cf.~\cite{zessin05}
or \cite{diparbeit}) subset of the product space
$\mathfrak{X}_d\times\mdrd$. The elements $(\x,\eta)\in\mathscr{D}$
are the clusters and point conf\/igurations which are `connected' in the
context of the connection rule $\mathscr{D}$.
Although in general misleading, it might be helpful to imagine the
(random) points $\eta$ as a set of nuclei, and a connected cluster $X$
as the orbiting electrons of one nucleus, where - of course - not all
conf\/igurations are possible.

If we take Voronoi and Delone tessellations as typical objects we want
to describe, it is easy to see that, for a given point conf\/iguration
$\eta\in\mdrd$, we need two concepts for the clusters. In the case of the
Voronoi tessellation the vertices of the cells are dif\/ferent from the
generating point conf\/iguration, while the vertices
of the cells in the Delone case coincide with the generating points. To give those
concepts a name we will call $\x\in\mathfrak{X}_d$ a \emph{cluster (of Type
  $\mathscr{D}$) \textbf{for} $\eta$}  if just
$(\x,\eta)\in\mathscr{D}$. If, additionally, $\x\subset\eta$ we will 
call $\x$ a \emph{cluster \textbf{in} $\eta$}. In this notation, the
vertices of the Voronoi cells are just clusters \emph{for} the
generating point conf\/igurations, while the vertices of the Delone
cells are clusters \emph{in} it. (To stay in the image of nuclei,
clusters \emph{in} a conf\/iguration might be imagined as
neighboured or even interacting nuclei, where the cluster properties
are the interaction rules.) 

It can also be convenient to extend the concept of cluster properties
to a situation where the clusters and the point conf\/iguration do not
lie in the same underlying space $\Rd$. Therefore we will also consider measurable
subsets of $\mathfrak{X}_d\times\mdrn$ as cluster properties, where
$d\not=n$. In this context, only clusters \emph{for} a point
conf\/iguration are well def\/ined. We will see an example below.

\begin{figure}
   \centering
   \includegraphics[width=\colwidth]{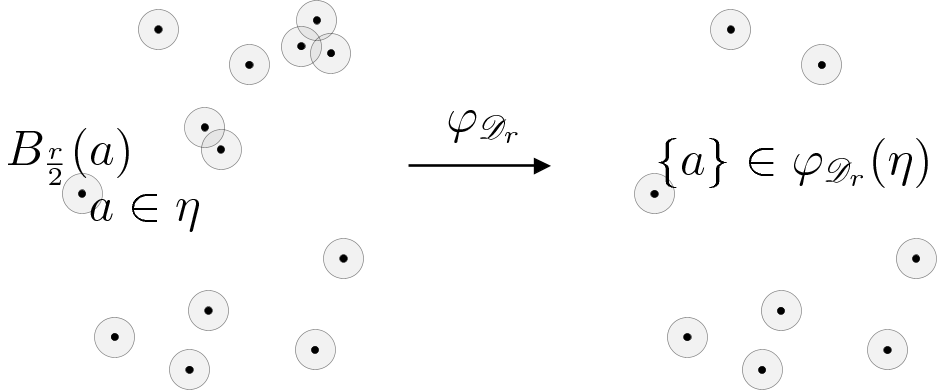}
   \caption{The clusters of type $\mathscr{D}_r$ in a point
     conf\/iguration 
     might be interpreted as a collection of non intersecting balls}
   \label{fig:hardcore}
\end{figure}

Before we get to more dif\/f\/icult examples and tilings, let us consider a
simple cluster property for illustration. Let
$r\in\R^+$ be f\/ixed. The pair $(\x,\eta)$ is def\/ined to belong to the cluster
property $\mathscr{D}_r\subset\mathfrak{X}_d\times\mdrd$ if $\x$
consists of a single point $a_\x$ and 
$\eta$ has no point which is closer to $a_\x$ than $r$ except possibly
$a_\x$, when this
is also an element of $\eta$.  
One cluster of type $\mathscr{D}_r$ \emph{for} a given conf\/iguration
$\eta$ is easily interpreted as a ball of radius $r$ which does not
intersect with the points of $\eta$. But the collection of clusters for
$\eta$ has no reasonable interpretation. The collection is also not
locally f\/inite in the sense that only f\/initely many clusters intersect
with arbitrary bounded sets of $\Rd$.
On the other hand, the collection of clusters  \emph{in} $\eta$ is
locally f\/inite and can be interpreted as a collection of balls of
radius $r/2$ around the points of $\eta$ where all the intersecting
balls are removed. Imagine $\eta$ now as a random realization through
some point process. Then, the collection of clusters \emph{in} $\eta$
`is' a random collection of non intersecting balls 
(see Figure \ref{fig:hardcore}). 
Note that in the general case, where clusters in or for a point
conf\/iguration consist of more than one point, the clusters may intersect.

If one wants to examine collections of clusters of certain type for or
in random point conf\/igurations, it is interesting to know whether the
underlying probability law, namely the chosen point process, produces
interesting collections of clusters. One indicator for this would be
inf\/initely (but locally f\/initely) many clusters for almost all
point conf\/igurations (with respect to the point process). In the
case of the clusters in a 
random conf\/iguration, there is the following interesting result
from \cite{zessin05}:  

\begin{zeroinfty}
  Suppose that the cluster property is of the kind that translating a
  pair $(\x,\eta)$ does not alter whether it belongs to the cluster
  property or not. Let also $P$ be a stationary point process. 
  Then, with probability one, we find either infinitely
  many clusters in $\eta$ or none.
\end{zeroinfty}

In the case of the Poisson point process $P_\lambda $, and subject to a mild
extra assumption (see \cite{zessin05}) one has an even stronger result: 
Assume that with positive probability there exists at least one
cluster in a randomly realized point conf\/iguration. This is suf\/f\/icient
to almost
always (with respect to $P_\lambda$) having inf\/initely many clusters
in such a conf\/iguration.
Again, it is convenient and possible to express collections, this time
not of points but of clusters, by means of Dirac measures. Here,
the collection $\varphi_{\mathscr{D}}(\eta)$ of clusters in
$\eta$ might be expressed by    
\begin{eqnarray}
  \label{eq:cluster-function}
   \varphi_\mathscr{D}(\eta)
  & = &\sum_{\x\text{ is a cluster of type $\D$ in }\eta}\delta_\x\notag\\
  & = &\sum_{\x\subseteq\eta}1_{\D}(\x,\eta)\,\delta_\x\,,
\end{eqnarray}
where $1_{\D}$ denotes the indicator function of the cluster property $\D$.

Although we do not have the $0$-$\infty$-law in the case of clusters
\emph{for} a conf\/iguration, we will see an example where the collection of
all the clusters for random $\eta$'s is locally f\/inite.
If we combine this function $\varphi_\mathscr{D}$ with a point process
$P$ (more
precisely, we take the image of the point process under the
transformation $\varphi_\mathscr{D}$), we have a probability measure
$\varphi_\D (P)$ on cluster
conf\/igurations (again, see \cite{diparbeit} for details). If we see one
cluster as a `point' in $\mathfrak{X}_d$, the notation $\mdg$ for the
(locally f\/inite) cluster conf\/igurations is 
sensible. A probability measure on this spaces is called \emph{cluster
  process}. 

We will now collect some information about a special family of clusters.


\section{Geometry}
\begin{figure}[htb]
  \centering
  \includegraphics[width=\smallwidth]{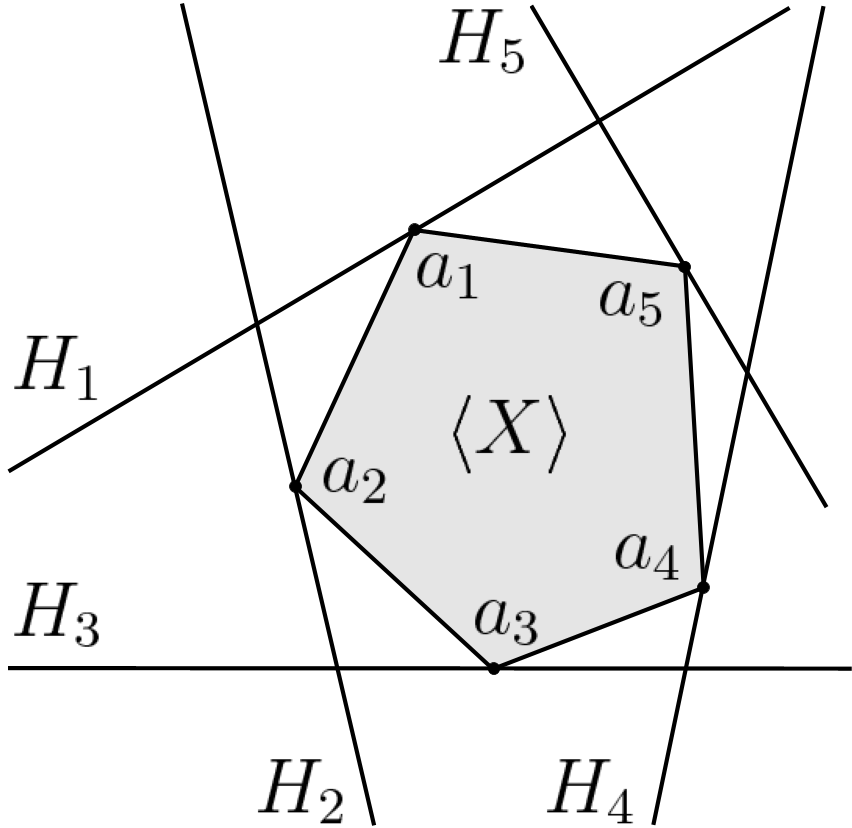}
  \caption{Illustration of a discrete polytope $\x=\{a_1,\ldots,a_5\}$
    and corresponding supporting hyperplanes $H_1,\ldots,H_5$} 
  \label{fig:disc-pol}
\end{figure}

The cluster conf\/igurations of interest for this article are
certain discretizations of tilings. Therefore, we need to adopt some well-known
concepts of tilings to this case. First we need to def\/ine
\emph{discrete polytopes}. Since we want to identify a convex polytope with
its vertex set, we take the properties of vertices for the
def\/inition: A cluster $\x\in\mathfrak{X}_d$ is a \emph{discrete polytope} if,
for all points $a\in \x$, there exists some hyperplane $H$ such that the
intersection of $H$ and the convex hull $\langle \x\rangle$ of $\x$
consists only of the point $a$. (See Figure \ref{fig:disc-pol} for
illustration.) 
It is quite obvious that there exists a one-to-one
correspondence between convex and discrete polytopes. The convex
polytope is retrieved from a discrete one by taking the convex
hull.

\begin{figure}
  \centering
  \label{fig:circumball}
  \includegraphics[width=\smallwidth]{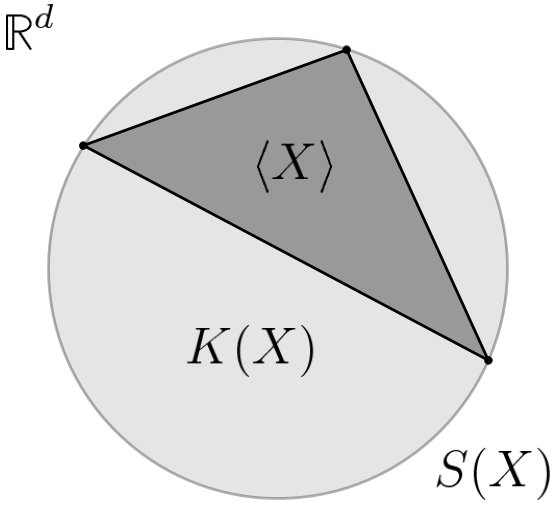}
  \caption{Circumball and circumsphere for a simplex (in this case a
    triangle)} 
\end{figure}

Similarly, \emph{discrete simplices} are obtained. For the
construction in the next section, we need a strong property of
simplices that is well-known in the $2$-dimensional case: Every
triangle has a uniquely def\/ined circumcircle. In  higher
dimensions, a full-dimensional simplex $\langle \x\rangle $ has a uniquely def\/ined
circumball $K(\x)$ where the complete vertex set -- the discrete
simplex $\x$ -- is contained in the border, the circumsphere $S(\x)$
(see Figure \ref{fig:circumball}); we refer to
\cite{diparbeit} for a proof. 

A collection $\mu\in\mdg$ of discrete polytopes is called a
\emph{cluster tessellation} if the collection of the convex hulls of
the clusters form a locally f\/inite \emph{face-to-face} tiling. 
In our context 'face-to-face' does not necessarily mean that there are
no holes in the tiling. It just states, that if two tiles intersect,
they intersect in whole faces.
The collection $\mu$ is
called \emph{simplicial} if all polytopes in it are simplices. It
is called \emph{complete} if the collection of convex hulls covers the
whole of $\Rd$. If a cluster process is concentrated on the set of
cluster tessellations, it is called \emph{random cluster tessellation}.


\section{Examples}
In this section, we give three examples for random cluster
tessellations, constructed via cluster properties. The f\/irst one
consists of clusters \emph{in} a point conf\/iguration and the second
one of clusters \emph{for} a point conf\/iguration. In both cases, the
underlying space for clusters and point conf\/igurations is the same,
the cluster properties are subsets of $\mathfrak{X}_d\times\mdrd$.
In the third example, the point conf\/igurations lie in some higher
dimensional space, where the underlying space of the clusters can be
interpreted as an embedded subspace.

\subsubsection*{First example: a special Delone tiling}
The cluster property that generates our tiling is def\/ined as follows:
Let $R\in\R^+$ be f\/ixed. For a discrete simplex, let
$\dot{K}(\x)=:K(\x)\smallsetminus \x$ denote the
circumball of  $\x$, where the discrete simplex is removed.

\begin{figure}[htb]
   \centering
   \includegraphics[width=\halfwidth]{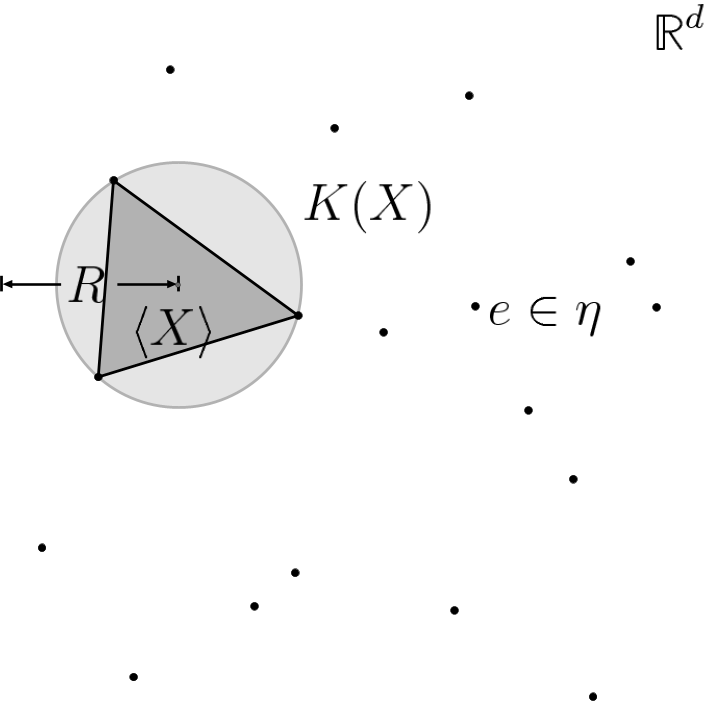}
   \caption{Illustration of clusters of type $\dr$ in a given $\eta$}
   \label{fig:del-cluster}
\end{figure}

\begin{figure}[htb]
   \centering
   \subfigure[$\x$ is not a cluster in $\eta$, because the radius of
   the circumball is to
   big]{\includegraphics[width=\halfwidth]{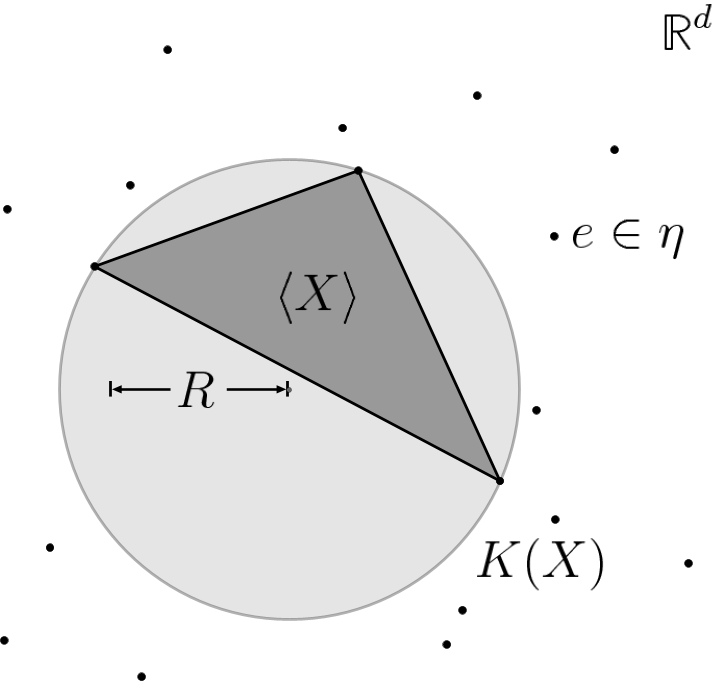}} 
   \hspace{0.5cm}
   \subfigure[$\x$ is a not cluster in $\eta$, because there is another
   point of the conf\/iguration in the
   circumball]{\includegraphics[width=\halfwidth]{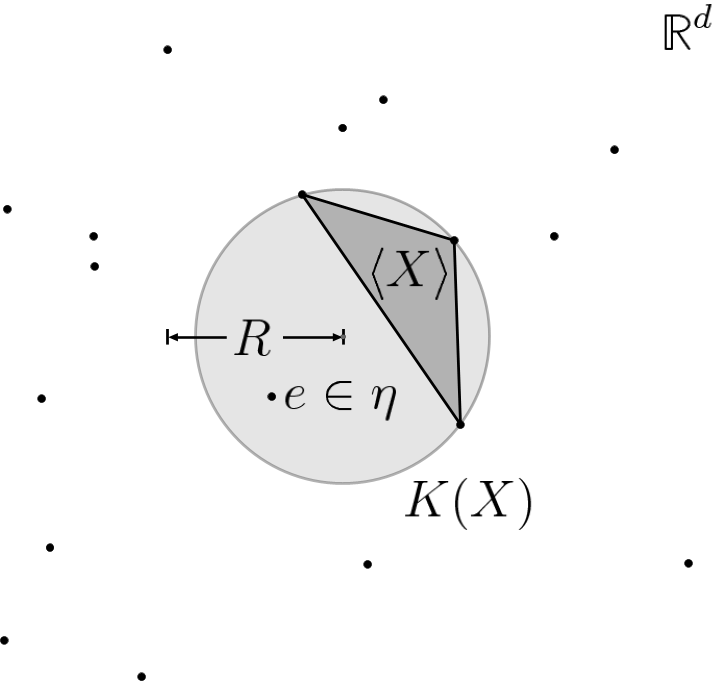}} 
   \caption{Illustration of clusters that fail to be of type $\dr$ in a
     given $\eta$} 
   \label{fig:no-del-cluster}
 \end{figure}

\begin{figure}[htb]
   \centering
   \includegraphics[width=\halfwidth]{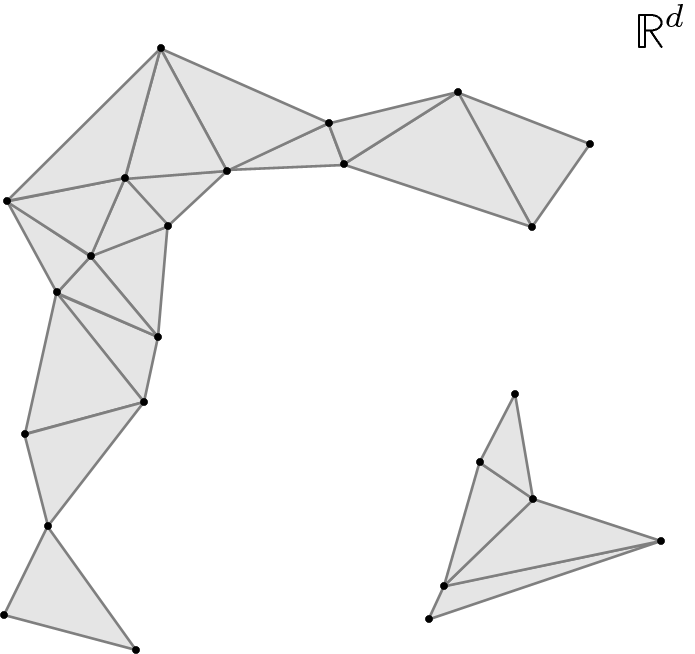}
   \caption{A typical section of a tessellation of the kind of Eq.~\eqref{eq:del-cluster}}
   \label{fig:delaunay}
\end{figure}

A tuple $(\x,\eta)\in\mathfrak{X}_d\times\mdrd$ belongs to the
cluster property $\dr$ if
\begin{enumerate}
 \item $\x$ is a $d$-dimensional simplex,
 \item $\dot{K}(\x)$ does not intersect with $\eta$, and
 \item $K(\x)$ has a radius $\leq R$.
\end{enumerate}
Figures \ref{fig:del-cluster} and \ref{fig:no-del-cluster} illustrate
clusters of type $\dr$ in a given point conf\/iguration $\eta$.
The f\/irst two assumptions are based on ideas of Delone
\cite{delaunay34} and ensure that, for a given $\eta\in\mdrd$, the
cluster conf\/iguration
\begin{equation}
  \label{eq:del-cluster}
  \varphi_{\dr}(\eta):=\sum_{\x\text{ is a cluster of type $\dr$ in }\eta}\delta_\x
\end{equation}
is face-to-face and simplicial. The assumption (iii) makes the
cluster conf\/iguration locally f\/inite and thus a tessellation. On the
other hand it produces holes in the tessellation (more precisely: in
the union of the convex hulls of the clusters) when the points in
$\eta$ are not dense enough.  (Figure \ref{fig:delaunay} gives a
typical section out of such an incomplete tessellation.)   
Thus, we get the following results:
\begin{proposition}
  Any simple point process $P$ generates a random tessellation in the form of
  $\varphi_{\dr} (P)$ as explained in Section \ref{sec:clusters}.
\end{proposition}
If we again consider the Poisson point process
$P_\lambda$, we can apply the $0$-$\infty$-law of stochastic geometry
to obtain:
\begin{proposition}
  $P_\lambda$-almost surely,
  \begin{enumerate}
  \item there are inf\/initely many clusters of type $\dr$ in a
    realization $\eta$, and
  \item there are holes in the tessellation.
  \end{enumerate}
\end{proposition}
Here, (ii) holds due to the fact, that 
a realization of a Poisson point process almost always has
arbitrarily big gaps somewhere between the points.
This, see \cite{diparbeit}, could
possibly lead to models for random holes in certain condensed matter. 

\subsubsection*{Second example: discrete Voronoi tilings}
This example describes a construction of the well-known Poisson
Voronoi tiling. In \cite{diparbeit}, a generalization, so-called
\emph{random Laguerre tessellations}, based on marked point processes
and the theory provided by Schlottmann \cite{schlottmann93},
is constructed. However, the idea of random tessellations constructed as
clusters \emph{for} random point conf\/igurations is better illustrated
in this easier case, so we stick to it. 

\begin{figure}[htb]
  \centering
  \includegraphics[width=\halfwidth]{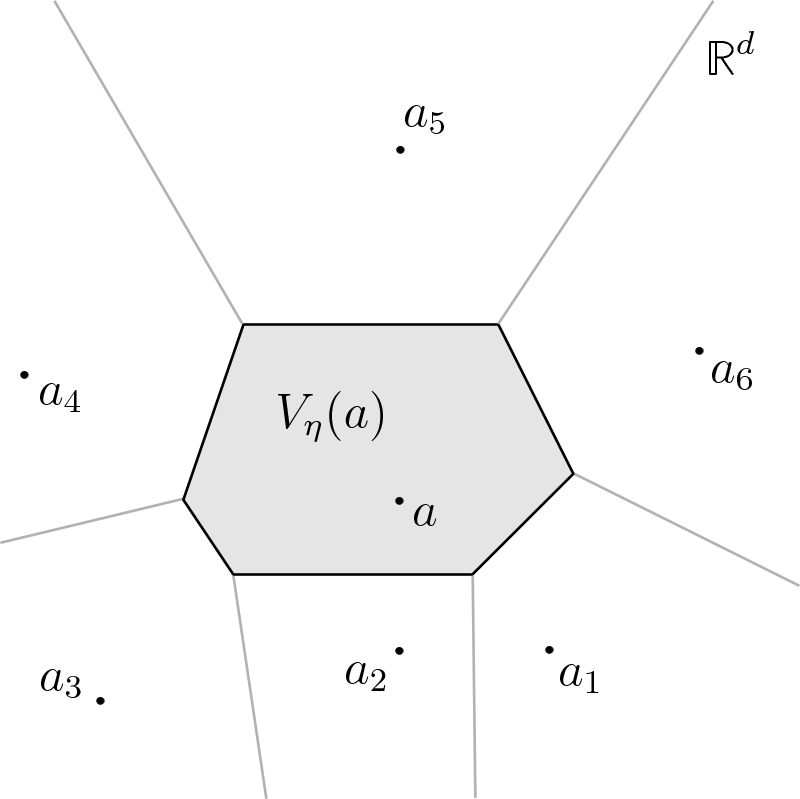}
  \caption{The Voronoi cell of $a$ in 
    $\eta=\delta_a+\delta_{a_1}+\delta_{a_2}+\ldots$}
  \label{fig:voronoi-cell}
\end{figure}

To understand the generating cluster property, recall the def\/inition of
Voronoi cells for a given point conf\/iguration $\eta\in\mdrd$. If
$a\in\eta$, the Voronoi cell of $a$ in $\eta$ is given by
\begin{equation*}
  V_\eta(a):=\{v\in\Rd\,\vert\,\text{No other
    point of }\eta\text{ is closer to }v\text{ than }a\}\,.
\end{equation*}
We call $a$ the \emph{center} of the Voronoi cell $V_\eta(a)$. 
The def\/inition is illustrated in Figure \ref{fig:voronoi-cell}.
If the convex hull of $\eta$ is the whole space $\Rd$, the Voronoi
cell is a convex polytope. 

\begin{figure}[htb]
  \centering
  \includegraphics[width=\halfwidth]{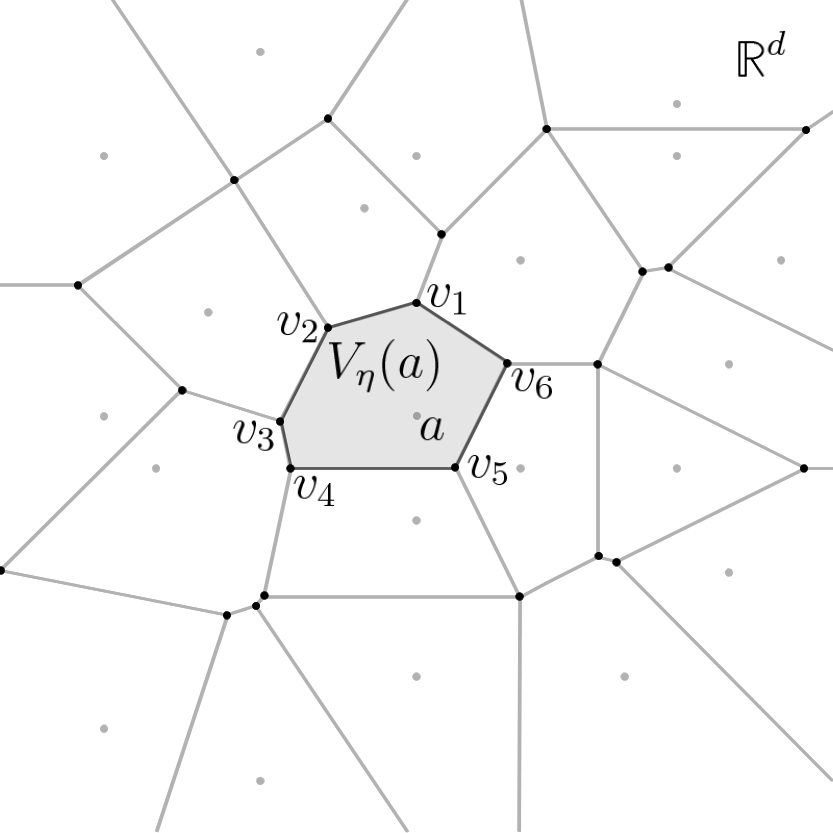}
  \caption{$\x=\{v_1,\ldots,v_6\}$ is a cluster of type $\dv$ for the
    conf\/iguration $\eta$ since it is the vertex set of the Voronoi
    cell $V_\eta(a)$} 
  \label{fig:voronoi2}
\end{figure}

We can now def\/ine our cluster property $\dv$ for this example:
$(\x,\eta)\in\dv $ if
\begin{enumerate}
\item the convex hull of $\eta$ is $\Rd$ and
\item $\x$ is the vertex set of a Voronoi cell in $\eta$.
\end{enumerate}
In this case, it is easy to see that, for every $\eta\in\mdrd$, there are
no clusters of type $\dv$ in $\eta$. But as long as (i) holds, the
collection 
\begin{equation*}
  \psi_{\dv}(\eta):=\sum_{\x\text{ is a cluster of type $\dv$ for }\eta}\delta_\x
\end{equation*}
is a complete cluster tessellation (a proof can be found in
\cite{schneider-weil00}). 
See Figure \ref{fig:voronoi2} for an illustration.
The Poisson point process $P_\lambda$ produces point conf\/igurations of
the kind (i) with probability one. Thus, we have:
\begin{proposition}
  $\psi_{\dv}(P_\lambda)$ is a complete random tessellation.
\end{proposition}

\subsubsection*{Third example: a random cut and project tiling}
This example is based on the so called \emph{cut and project scheme}, a
method to obtain tilings from a (higher dimensional) lattice.   
The vertices or the centers, in the sense of Voronoi cells, of the
tiles are projections of subsets of the lattice. 
A detailed description of the underlying theory can be found in \cite{moody97}. 
A large group of deterministic tilings can be constructed this way, for instance the
\emph{Penrose Tiling} (cf.~\cite{debruijn81}) or the
\emph{Ammann-Beenker Tiling} (cf.~\cite{beenker82}). We present a way
to construct random tilings which are `close' to the known 
deterministic ones, where `close' will have two meanings: in
the f\/irst one the random points still lie on the lattice but not every
point of the lattice will appear. The second interpretation will
produce one point for every lattice point, but the points might be
randomly shifted within some given radius.  

Again, we stick to a simple example, a  $1$-dimensional tiling. The
mechanisms for randomness can easily be adapted to any tiling
that can be obtained via the cut and
project scheme. The deterministic case of our example is
described in \cite{bagrmo02}. 

\begin{figure}[htb]
  \centering
  \includegraphics[width=\colwidth]{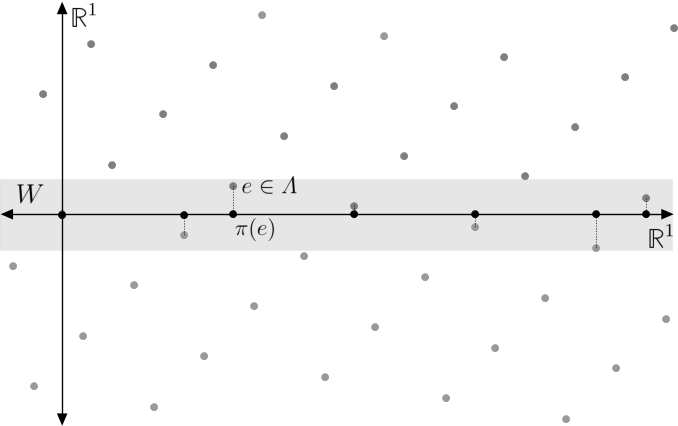}
  \caption{The projections of the points in $\varLambda\cap W$ form
    the vertex set of an aperiodic tiling} 
  \label{fig:silvermeans1}
\end{figure}

Consider the $2$-dimensional lattice
\begin{equation*}
  \varLambda:=\left\{\left(u+v\,\sqrt{2}\,,\,u-v\,\sqrt{2}\right)\middle\vert\,u,v
    \text{ integers}\right\}\subset\R^2\,.
\end{equation*}
Let $\pi$ be the projection onto the f\/irst coordinate and $\pi^*$ the
one onto the second. Consider the set
\begin{equation*}
  W:=\left\{e\in\R^2\,\middle\vert\,
    \pi^*(e)\in\left[-\frac{1}{\sqrt{2}},\frac{1}{\sqrt{2}}\right]\right\}\,, 
\end{equation*}
which we will call \emph{strip}.
The projection $\pi(\varLambda\cap W)$ forms the vertex set of the so-called
\emph{silver-mean chain}, which is a
deterministic 
\emph{aperiodic tiling} (again cf.\cite{bagrmo02}). Figure
\ref{fig:silvermeans1} illustrates this concept.

The appropriate cluster property $\ds\subset\mathfrak{X}_1\times\mdrt$
can be def\/ined as follows: 
$\left(X,\eta\right)\in\ds$ if and only if
\begin{enumerate}
\item $\card X=2$,
\item $X=\{a=\pi(e_1),b=\pi(e_2)\}$, with $e_1\not=e_2\in\eta\cap W$ and
\item the intersection of the open intervall
  $(a,b)$ and $\pi(\eta\cap W)$ is empty.
\end{enumerate}
In this $2$-dimensional case, it is easy to see that, as long as
$\pi(\eta\cap W)$ is a discrete point set, 
\begin{equation*}
  \psi_{\ds}(\eta):=\sum_{\x\text{ is a cluster of type $\ds$ for }\eta}\delta_\x
\end{equation*}
is a $1$-dimensional tessellation where neighboured points form a
cluster, respectively the vertex set of a tile.

To embed the deterministic version of the tiling into point process
theory, just def\/ine $P_\varLambda$ to be the point process in $\R^2$ which
produces the lattice $\varLambda$ with probability one. Then,
$\psi_{\ds}(P_\varLambda)$ is a process that almost surely produces
the silver means tiling. This tiling consists of two prototiles of
length $1$, and $1+\sqrt{2}$, respectively. 

\begin{figure}[htb]
  \centering
  \includegraphics[width=\colwidth]{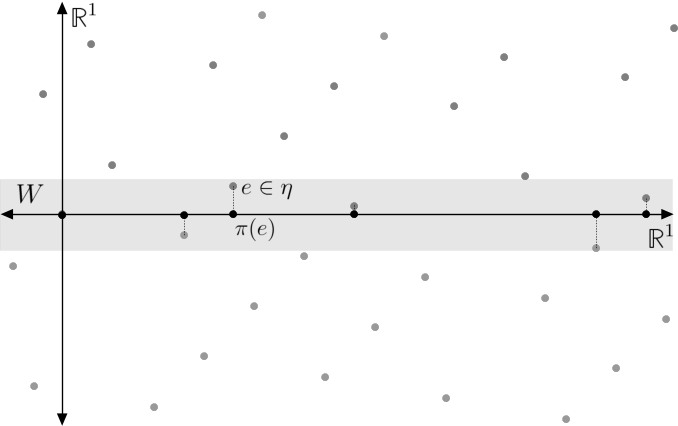}
  \caption{$P_\rho^*$ produces random subsets of the lattice
    $\varLambda$, the tiles might become larger}
  \label{fig:silvermeans3}
\end{figure}

For the f\/irst random version, consider the discrete measure
\begin{equation*}
  \rho:=\sum_{e\in\varLambda}c\cdot\delta_e\,,
\end{equation*}
where $c$ is some positive constant.
As mentioned above, $P_\rho$, the Poisson point process with intensity
measure $\rho$, might produce point conf\/igurations with more than one
atom in a single point, in this case in the points of
$\varLambda$. All the point sets produced by the random mechanism $P_\rho$ have
the form
\begin{equation*}
  \eta=\sum_{e\in\varLambda}k_e(\eta)\cdot\delta_e\,,
\end{equation*}
where $k_e(\eta)$ is some natural number or $0$. The \emph{support} of
such an
$\eta$ is def\/ined as
\begin{equation*}
  \supp (\eta):=\sum_{e\in\varLambda,\,k_e(\eta)\not=0}\delta_e\,.
\end{equation*}
The support is a subset of $\varLambda$, especially has only one atom
at every point.  
Let $P^*_\rho $ be the image of $P_\rho $ under the mapping $\supp $. Thus:
\begin{proposition}
  $P^*_\rho$ is again a simple point process, where the realizations
  are random subsets of the lattice. The probability for a certain
  point $e$ of $\varLambda$ to be in the random set is $1-\e^{-c}$.
\end{proposition}
Figure \ref{fig:silvermeans3} shows a typical
randomly realized point set and what happens by taking the cut and project clusters. 
Since the holes in such a random subset cannot be controlled, the
tiles corresponding to the clusters for a random $\eta$ might have any
length of the form $n+m\sqrt{2}$, $n,m$ non-negative integers.

The next point process randomly shifts the points of the lattice. Another way to randomly shift the points is presented in
\cite{hof95}. There might be some way to transform these two
approaches into one another.
For the construction of our point process, let $\varepsilon>0$ and
$B_\varepsilon(e)$ 
the ball of radius $\varepsilon$ and centre $e$,
$e\in\varLambda$. $\varepsilon$ should be choosen small enough so that the
balls do not intersect. Consider the mapping
$\boldsymbol{b}_e:\mdrt\to\mdrt$,
\begin{eqnarray*}
  \boldsymbol{b}_e(\eta):=
  \begin{cases}
    \frac{1}{\card(\eta\cap B_\varepsilon(e))}\sum_{f\in\eta\cap
      B_\varepsilon(e)}f\,, & \text{if }\eta\cap
    B_\varepsilon(e)\not=\emptyset\,,\\
    e\,,& \text{if }\eta\cap B_\varepsilon(e)=\emptyset\,,
  \end{cases}
\end{eqnarray*}
which gives the barycentres of all the points in $\eta\cap
B_\varepsilon(e)$.
The conf\/iguration
\begin{equation*}
  \boldsymbol{b}(\eta):=\sum_{e\in\varLambda}\boldsymbol{b}_e(\eta)
\end{equation*}
of all these barycentres, by construction, has exactly one point in every
$\varepsilon$-ball around the lattice points.
\begin{figure}[htb]
  \centering
  \includegraphics[width=\colwidth]{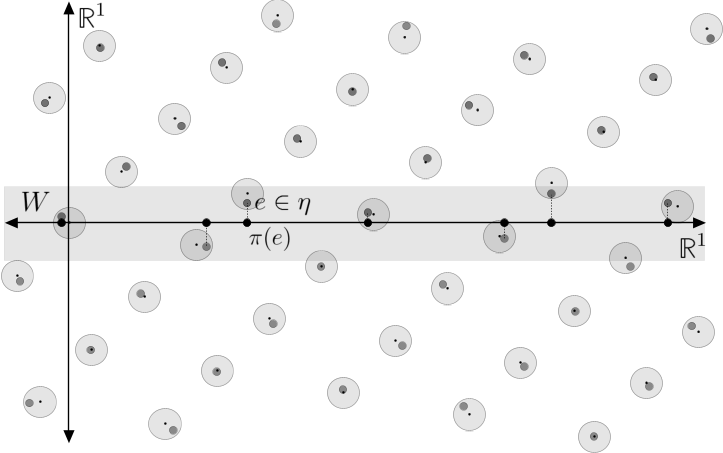}
  \caption{For typical realisations $\eta$ of $P_{\boldsymbol{b}}$ the
    projections of $\eta\cap W$ form slightly changed tiles or even
    completely new ones} 
  \label{fig:silvermeans2}
\end{figure}
Thus we have the following result:
\begin{proposition}
  If you take some arbitrary simple point process $P$, e.g.
  $P_\lambda$, the image $P_{\boldsymbol{b}}$ of $P$ under the mapping
  $\boldsymbol{b}$ randomly produces point conf\/igurations with
  exactly one point in every $\varepsilon$-ball centred in the lattice
  points.
\end{proposition}
Figure
\ref{fig:silvermeans2} illustrates the typical situation and the
resulting projections. The corresponding tiles to the clusters for a given
$\eta$ typically are close to the tiles of the original silver means
tiling, differing in length up to $2\cdot\varepsilon$. But since
certain $\varepsilon$-balls of points in $\varLambda\cap W^\complement $
intersect with the strip $W$, there might be `completely new'
tiles. Nevertheless, the density of the vertices stays the same, since
the probabilities to shift a point into and to the outside of the strip are the same.

In both of the cases of the third example of this article, it is easy
to see the following:
\begin{proposition}
  The cluster processes $\psi_{\ds}(P^*_\rho)$, respectively
  $\psi_{\ds}(P_{\boldsymbol{b}})$, are random $1$-di\-men\-sional
  tessellations.
\end{proposition}


\section{Conclusions}
Point processes, especially the Poisson point processes, in
combination with cluster properties give access to modelling discrete
random structures. The information of the cluster properties in this
case carry the information about possible connections respectively
interactions of the particles. The third type of example
shows a way to slightly randomise aperiodic tilings. In this context, future
applications to random tilings in the sense of Gummelt \cite{gummelt04,gummelt06}
are of interest. If point processes could be constructed that almost
surely produce more special configurations, e.g. Delone or FLC
(cf. for instance \cite{lms02}) sets, the
presented methods might get closer to applications like glasses or foams.


\subsubsection*{Acknowledgements}
It is a pleasure to thank M.~Baake, D.~Frettlöh, C.~Richard and
H.~Zessin for several useful comments and clarifying discussions. I
also want to thank the reviewers for a number of very helpful
suggestions to improve the manuscript.

This work was partially supported by the German Research
Council (DFG), within the CRC 701.


\end{document}